\documentclass{easychair}

\usepackage{lineno}
\usepackage{doc}
\usepackage{booktabs} 
\usepackage{amsmath,amsfonts}
\usepackage{float}
\usepackage{algorithm}
\usepackage{algpseudocode}
\usepackage{graphicx}
\usepackage{multirow}
\usepackage{multicol}
\usepackage{amssymb}
\usepackage{subcaption} 
\usepackage{xcolor}
\usepackage{diagbox}
\usepackage{enumitem}
\setlist[itemize]{leftmargin=*}
\usepackage{url}
\usepackage{threeparttable}
\usepackage{pifont}
\usepackage{nopageno}
\usepackage{bm}
\usepackage{scalerel}
\usepackage{makecell}
\usepackage[colorinlistoftodos]{todonotes}

\newcommand{\R}{\mathbb{R}}

\title{MORCIC: Model Order Reduction Techniques for Electromagnetic Models of Integrated Circuits}

\author{
Dimitrios Garyfallou\inst{1},
Athanasios Stefanou\inst{2},
Christos Giamouzis\inst{1},
Moschos Antoniadis\inst{2},
Georgios Chararas\inst{2},
Konstantinos Chatzis\inst{2},
Dimitris Samaras\inst{2},
Rafaela Themeli\inst{2},
Anastasios Michailidis\inst{2},
Vasiliki Gogolou\inst{2},
Nikos Zachos\inst{2},
Nestor Evmorfopoulos\inst{1},
Thomas Noulis\inst{2},
Vasilis F. Pavlidis\inst{2},
Alkiviadis Hatzopoulos\inst{2},
Elpida Chatzineofytou\inst{3},
\and
Yiannis~Moisiadis\inst{3}
}

\institute{
 Dept. of Electrical and Computer Engineering, 
University of Thessaly, Volos, Greece \\
  \email{\{digaryfa, cgiamouzis, nestevmo\}@e-ce.uth.gr}
\and
   Aristotle University of Thessaloniki, Thessaloniki, Greece\\
   \email{\{tnoul, vpavlid, alkis\}@auth.gr}\\
\and
  ANSYS-Hellas, Athens, Greece\\
   \email{\{elpida.chatzineofytou,  yiannis.moisiadis\}@ansys.com}\\
}

\begin{document}
\maketitle

\begin{abstract}
Model order reduction (MOR) is 
crucial
for the design process of integrated circuits. 
Specifically, the vast amount of passive RLCk elements in electromagnetic models extracted from physical layouts exacerbates the extraction time, the storage requirements, and, most critically,~the post-layout simulation time of the analyzed circuits. 
The MORCIC project aims to overcome 
this problem by proposing new MOR techniques that perform better 
than commercial tools. 
Experimental evaluation on several analog and mixed-signal circuits with millions of 
elements indicates that the proposed methods lead to $\times$5.5 smaller ROMs 
while maintaining 
similar accuracy compared to
golden ROMs provided by ANSYS~RaptorX\texttrademark.
\end{abstract}
\section{Introduction}
\label{sec:intro}
Electromagnetic model extraction plays a key role in the design and analysis of integrated circuits. The extracted models are simulated to accurately predict the behavior of the passive elements of the design. Model order reduction (MOR) can reduce the complexity of RLCk models (consisting of resistances, capacitances, self-inductances, and mutual inductances)
with many elements ($>$1M) and ports ($>$10), while retaining an accurate approximation of the input and output behavior of the circuit~\cite{Odabasioglu1998, Antoulas2004}. 
In this way, the simulation time of complex systems can be decreased to a few seconds by constructing systems of smaller dimensions that preserve the essential characteristics of the original model.

MOR methods are distinguished into two main categories. Moment matching (MM) techniques~\cite{Odabasioglu1998} are typically preferred due to their computational efficiency. 
However, they rely on an ad-hoc selection of the number of moments to obtain accurate reduced-order models (ROMs). Moreover, the ROM order is highly dependent on the number of ports, posing significant challenges to the reduction of multi-port models.
Contrary to MM methods,~system theoretic techniques, like balanced truncation (BT)~\cite{Antoulas2004}, offer reliable error bounds and have no fundamental limitation to the number of ports they can handle.
Nevertheless, BT applies only to small-scale models since it involves the solution of computationally expensive Lyapunov~equations~\cite{Antoulas2004}. 

During the MORCIC project, appropriate performance improvements are explored to overcome the main drawback of the conventional BT method. 
To this end, we adopt an efficient low-rank technique based on the extended Krylov subspace (EKS) for solving the Lyapunov equations.
The proposed approach can be integrated into industrial extraction tools, such as the ANSYS RaptorX\texttrademark ~\cite{raptX}, to obtain more compact ROMs of large-scale multi-port RLCk models.
\section{Background}
\label{sec:background}
An RLCk circuit with $n$ nodes, $m$ inductive branches, $p$ inputs, and $q$ outputs can be described in the time domain via the Modified Nodal Analysis (MNA) as~\cite{MNA}: 
\begin{equation} \label{mna}
    \begin{aligned}
        \begin{pmatrix}
            \mathbf{G_{n}} & \mathbf{E} \\
            \mathbf{-E}^T & \mathbf{0} 
        \end{pmatrix}
        \begin{pmatrix}
            \mathbf{v}(t) \\
            \mathbf{i}(t) 
        \end{pmatrix} +
        \begin{pmatrix}
            \mathbf{C_{n}} & \mathbf{0} \\
            \mathbf{0} & \mathbf{M} 
        \end{pmatrix}
        \begin{pmatrix}
            \dot{\mathbf{v}}(t) \\
            \dot{\mathbf{i}}(t) 
        \end{pmatrix} = 
        \begin{pmatrix}
            \mathbf{B}_1 \\
            \mathbf{0} 
        \end{pmatrix}
        \mathbf{u}(t), \quad 
        \mathbf{y}(t) = 
        \begin{pmatrix}
            \mathbf{L}_1 \quad
            \mathbf{0} 
        \end{pmatrix}
        \begin{pmatrix}
            \mathbf{v}(t) \\
            \mathbf{i}(t) 
        \end{pmatrix}
    \end{aligned}
\end{equation}
where $\mathbf{G_{n}} \in\R^{n\times n}$ is the node conductance matrix,  $\mathbf{C_{n}}\in\R^{n\times n}$ is the node capacitance matrix, $\mathbf{M}\in\R^{m\times m}$ is the branch inductance matrix, $\mathbf{E}\in\R^{n \times m}$ is the node-to-branch incidence matrix, $\mathbf{v}\in\R^{n}$ is the vector of node voltages, $\mathbf{i}\in\R^{m}$  is the vector of inductive branch currents, $\mathbf{u} \in \R^{p}$  is the vector of input excitations, $\mathbf{B}_1\in\R^{n\times p}$ is the input-to-node connectivity matrix, $\mathbf{y}\in\R^{q}$ is the vector of output measurements, and $\mathbf{L}_1\in\R^{q\times n}$ is the node-to-output connectivity matrix. 
Moreover, we denote
$\dot{\mathbf{v}}(t) \equiv \frac{d\mathbf{v}(t)}{dt}$ and $\dot{\mathbf{i}}(t) \equiv \frac{d\mathbf{i}(t)}{dt}$.

\noindent
If we now define the model order as $N \equiv n + m$,
the state vector as $\mathbf{x}(t) \equiv \begin{pmatrix}
\mathbf{v}(t) \\
\mathbf{i}(t) 
\end{pmatrix}$, and~also:
\begin{equation*}
\begin{aligned}
\mathbf{G}\equiv -\begin{pmatrix}
\mathbf{G_{n}} & \mathbf{E} \\
\mathbf{-E}^T & \mathbf{0} 
\end{pmatrix},\quad  \mathbf{C} \equiv \begin{pmatrix}
\mathbf{C_{n}} & \mathbf{0} \\
\mathbf{0} & \mathbf{M} 
\end{pmatrix},\quad \mathbf{B}\equiv \begin{pmatrix}
\mathbf{B}_1 \\
\mathbf{0} 
\end{pmatrix}, \quad \mathbf{L}\equiv \begin{pmatrix}
\mathbf{L}_1 \quad
\mathbf{0} 
\end{pmatrix}
\end{aligned}, 
\end{equation*}
then Eq.~(\ref{mna}) can be written in the 
state-space form of an $N$-th order linear dynamical system:
\begin{equation}
\begin{aligned} \label{state}
\mathbf{C}\frac{d \mathbf{x}(t)}{d t} = \mathbf{G x}(t) + \mathbf{Bu}(t), \quad 
\mathbf{y}(t) = \mathbf{L x}(t).
\end{aligned}
\end{equation}
\noindent The objective of MOR is to produce an approximate ROM: 
\begin{equation}
\begin{aligned}
\mathbf{\tilde C} \frac{d \mathbf{\tilde x}(t)}{d t} =\mathbf{\tilde G} \mathbf{\tilde x}(t) + \mathbf{\tilde B} \mathbf{u(t)}, \quad
\mathbf{\tilde y}(t) = \mathbf{\tilde L \tilde x}(t)
\end{aligned}
\end{equation}				
where  $\mathbf{\tilde G}, \mathbf{\tilde C} \in \R^{r\times r} $, $\mathbf{\tilde B} \in \R^{r\times p} $, $\mathbf{\tilde L} \in \R^{q\times r}$, the reduced order $r<<N$, and the ROM approximates the original model in the sense that the output error is bounded as $||\mathbf{\tilde{y} }(t) -\mathbf{y}(t)||_2 < \varepsilon||\mathbf{u}(t)||_2$ for given  $\mathbf{u}(t)$ and small threshold $\varepsilon$. The output error bound can be 
equivalently expressed in the frequency domain as $||\mathbf{\tilde{y} }(s) -\mathbf{y}(s)||_2 < \varepsilon||\mathbf{u}(s)||_2$ 
(via 
Plancherel's theorem~\cite{Plancherel}).~If
\begin{equation*}
\begin{aligned}
\mathbf{H}(s) = \mathbf{L}(s\mathbf{C} - \mathbf{G})^{-1} \mathbf{B}, \quad
\mathbf{\tilde H}(s) =  \mathbf{\tilde L}(s\mathbf{\tilde C} - \mathbf{\tilde G})^{-1} \mathbf{\tilde B}
\end{aligned}
\end{equation*}
are the transfer functions of the original model and the ROM, the 
output
error~becomes: 
\begin{equation*}
\begin{aligned}
||\mathbf{\tilde{y} }(s) -\mathbf{y}(s)||_2 = ||\mathbf{\tilde{H}}(s) \mathbf{u}(s) - \mathbf{H}(s)\mathbf{u}(s)||_2 \quad \leq \quad ||\mathbf{\tilde{H}}(s) - \mathbf{H}(s)||_\infty||\mathbf{u}(s)||_2	
\end{aligned}
\end{equation*}
where $||.||_\infty$ is the 
$\mathcal{L}_2$ matrix norm (or $\mathcal{H}_\infty$ norm of a rational transfer function). 
Thus,~to~bound this error,
we must bound the distance between the transfer functions as:~$||\mathbf{\tilde{H}}(s)~-~\mathbf{H}(s)||_\infty~<~\varepsilon$.

\noindent
To achieve this, BT transforms the original model into a ROM with a "balanced" state vector and then truncates the joint controllability-observability singular values of the system (so-called Hankel singular values) that sum up to the given threshold $\varepsilon$, as described in~\cite{Antoulas2004}.

\section{MOR by Balanced Truncation}
\label{sec:bt}
\subsection{Initial BT MOR}
\label{sec:bt_vanila}
BT  relies on the computation of the controllability Gramian $\mathbf{P}$ and the observability Gramian~$\mathbf{Q}$,
which are calculated as the solutions of the following Lyapunov matrix equations \cite{Antoulas2004}:
\begin{equation}
\begin{aligned}\label{Eq:lyap_sta}
(\mathbf{C}^{-1}\mathbf{G}) \mathbf{P} +  \mathbf{P} (\mathbf{C}^{-1}\mathbf{G})^T  = - (\mathbf{C}^{-1}\mathbf{B}) (\mathbf{C^{-1}}\mathbf{B})^T, \quad
(\mathbf{C}^{-1}\mathbf{G})^T \mathbf{Q} +  \mathbf{Q}(\mathbf{C}^{-1}\mathbf{G}) = - \mathbf{L}^T \mathbf{L}.
\end{aligned}
\end{equation}
The main steps of the BT procedure are summarized in Algorithm \ref{vanilla_BT}. 
As can be seen, the operations involved (e.g., the solution of Lyapunov equations and the singular value decomposition [SVD]) are computationally expensive with complexity $O(N^3)$. Moreover, they are applied on dense matrices, since the Gramians $\mathbf{P}, \mathbf{Q}$ are dense even if the system matrices $\mathbf{C}, \mathbf{G}, \mathbf{B}, \mathbf{L}$ are sparse. Consequently, 
the significant computational and memory cost for deriving the ROM hinders the applicability of BT to large-scale models (with order $N$ over a few thousand states).
\begin{figure}[!hbt]
\vspace{-2em}
\begin{algorithm}[H]
\small
	\caption{MOR by balanced truncation}\label{vanilla_BT}
        \textbf{Input:}  $ \mathbf{G}, \mathbf{C}, \mathbf{B}, \mathbf{L} $\\
	\textbf{Output:} $\mathbf{\tilde G}, \mathbf{\tilde C}, \mathbf{\tilde B}, \mathbf{\tilde L}$
	\begin{algorithmic}[1]
		\State Solve the Lyapunov equations to obtain the Gramian matrices $\mathbf{P}$ and $\mathbf{Q}$~\cite{Lathauwer2004}
		\State Compute the
            SVD of the Gramian matrices: $\mathbf{P} = \mathbf{U}_P \mathbf{\Sigma}_P \mathbf{V}_P^T $ and $\mathbf{Q} = \mathbf{U}_Q \mathbf{\Sigma}_Q \mathbf{V}_Q^T $
            \State Find the square root of the Gramian matrices: $\mathbf{Z}_P = \mathbf{U}_P \mathbf{\Sigma}_P^{1/2}$ and $\mathbf{Z}_Q = \mathbf{U}_Q \mathbf{\Sigma}_Q^{1/2}$
            \State Compute the SVD of the product of the roots: $\mathbf{Z}_Q^T\mathbf{Z}_P = \mathbf{U}\mathbf{\Sigma}\mathbf{V}^T$	
  \State  
  Compute transformation matrices: $\mathbf{T}_{(r\times N)}$\ =\ $\mathbf{\Sigma}_{(r\times r)}^{-1/2} \mathbf{U}_{(r\times N)} \mathbf{Z}_Q^T $,\ $\mathbf{T}_{(N\times r)}^{-1}$\ =\ $ \mathbf{Z}_P\mathbf{V}_{(N\times r)}\mathbf{\Sigma}^{-1/2}_{(r\times r)}$
            \State Compute 
            ROM: $\mathbf{\tilde G}$\ =\ $\mathbf{T}_{(r\times N)}\mathbf{G}\mathbf{T}_{(N\times r)}^{-1}$, \ $\mathbf{\tilde C}$\ =\ $\mathbf{T}_{(r\times N)}\mathbf{C}\mathbf{T}_{(N\times r)}^{-1}$,\  $\mathbf{\tilde B}$\ =\ $\mathbf{T}_{(r\times N)}\mathbf{B},\ \mathbf{\tilde L}$\ =\ $\mathbf{L}\mathbf{T}_{(N\times r)}^{-1}$
	\end{algorithmic}
\end{algorithm}
\vspace{-1.5em}
\end{figure}

\noindent
However, the products $(\mathbf{C}^{-1}\mathbf{B}) (\mathbf{C^{-1}}\mathbf{B})^T$ and $\mathbf{L}^T\mathbf{L}$ 
have low numerical order compared~to~$N$, as
$p,q<<N$, resulting in low-rank Gramian matrices that can be approximated using low-rank techniques. This greatly reduces the complexity and memory requirements of the solution of the Lyapunov equations and the SVD analysis, which are now of order $k$ instead of full order~$N$.

\subsection{Low-rank BT MOR}
\label{sec:bt_eks}
The essence of low-rank BT MOR is to iteratively project the Lyapunov Eq.~(\ref{Eq:lyap_sta}) onto a lower-dimensional Krylov subspace ($\mathcal{K}_k$)~\cite{Simoncini2007} and then solve the resulting small-scale equations to obtain low-rank approximate solutions of Eq.~(\ref{Eq:lyap_sta}). 
In this work, we exploit the EKS to accelerate the convergence to the final solution~\cite{ASPDAC21}.
The complete EKS method is presented in Algorithm~\ref{eksm-alg}. 
\begin{figure}[!hbt]
\vspace{-2em}
\begin{algorithm}[H]
\small
\caption{Extended Krylov subspace method for low-rank solution of Lyapunov equations}\label{eksm-alg}	
	\textbf{Input:}  $ \mathbf{G}_{C} \equiv \mathbf{C}^{-1}\mathbf{G},  \mathbf{B}_{C} \equiv \mathbf{C}^{-1}\mathbf{B}$  (or  $\mathbf{G}_{C}^T$, $\mathbf{L}^T $)\\
	\textbf{Output:} $\mathbf{Z}$ such that $\mathbf{P} \approx \mathbf{Z} \mathbf{Z} ^T $
	\begin{algorithmic} [1]
		\State $j=1$; $p=size\_col(\mathbf{B}_{C})$
        \label{eksm-alg:step2}
		\State $\mathbf{K}^{(j)} = Orth([\mathbf{B}_{C},\mathbf{G}_{C}^{-1}\mathbf{B}_{C}])$
		\While {$j < maxiter$}
        \label{eksm-alg:step4}
    		\State $\mathbf{A} = \mathbf{K}^{(j)T}\mathbf{G}_{C}\mathbf{K}^{(j)} $;\quad $\mathbf{R} = \mathbf{K}^{(j)T}\mathbf{B}_{C} $  
        \label{eksm-alg:step5}
    		\State Solve $\mathbf{A}\mathbf{X} +\mathbf{X} \mathbf{A}^{T}  = -\mathbf{R}\mathbf{R}^{T}$ for  $\mathbf{X} \in \R^{2pj \times 2pj}$
    		\If {converged}
        		\State $[\mathbf{U}, \mathbf{\Sigma}, \mathbf{V}] = \mathbf{SVD}(\mathbf{X})$; \quad
        		$\mathbf{Z} = \mathbf{K}^{(j)}\mathbf{U}\mathbf{\Sigma}^{1/2}$
        	    \State \textbf{break}
        	\EndIf
    		\State $k_1=2p(j-1)$; $k_2 = k_1+p$; $k_3 = 2pj$ 
        \label{eksm-alg:step11}
    		\State $\mathbf{K}_1  = [\mathbf{G}_{C}\mathbf{K}^{(j)}(:,k_1+1:k_2),\mathbf{G}_{C}^{-1}\mathbf{K}^{(j)}(:,k_2+1:k_3)]$
    		\State $\mathbf{K}_2 = Orth(\mathbf{K}_1) $ \quad w.r.t. \quad $\mathbf{K}^{(j)}$
    		\State $\mathbf{K}_3 = Orth(\mathbf{K}_2) $ 
    		\State $\mathbf{K}^{(j+1)} = [\mathbf{K}^{(j)},\mathbf{K}_3]$\quad 
		    \State $j=j+1$
		\EndWhile
	\end{algorithmic}
\end{algorithm}
\vspace{-2.5em}
\end{figure}

\section{Experimental Evaluation}
\label{sec:experimental}
\label{sec:exp_setup}
\begin{center}
\begin{figure}[!hbt]
\vspace{-0.5em}
  \centering
  \includegraphics[width=0.92\textwidth]{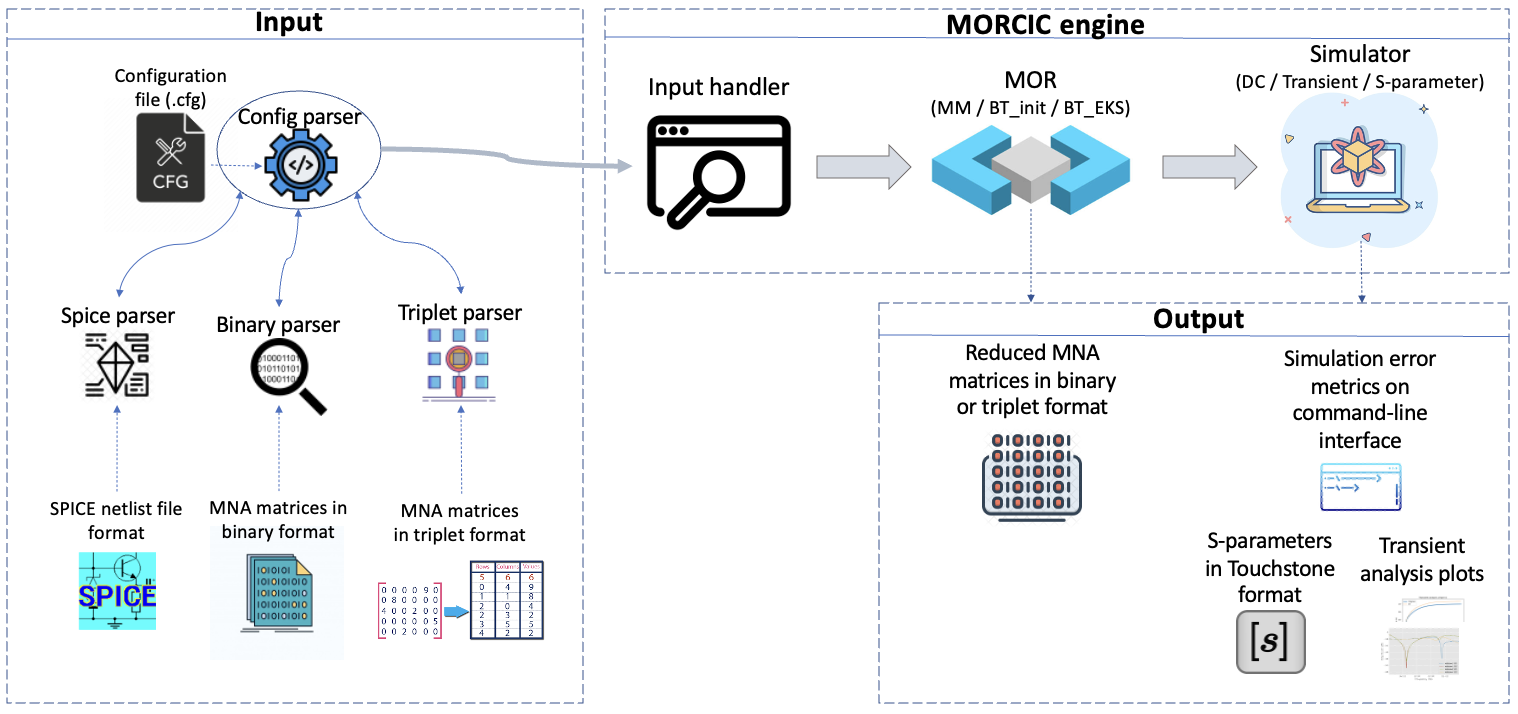}
  \vspace{-7pt}
  \caption{Software architecture of the MORCIC tool.}
  \label{fig:morcic_architecture}
\vspace{-2em}
\end{figure}
\end{center}

To evaluate the proposed MOR methods, we developed an EDA tool that implements 
the
BT algorithms 
presented in
Section~\ref{sec:bt}. As depicted in Figure~\ref{fig:morcic_architecture}, the only input is a configuration file that defines the path to the MNA matrices along with some parameters. After applying BT MOR, the tool performs DC, transient, and SP analysis, to compare the ROM to~the original model. The output 
includes the S-parameters and MNA matrices of the ROM. 
The cross-platform MORCIC tool was developed in C++ using the CMake automation software.
All experiments were executed on a Linux workstation with a 
2.80 GHz 16-thread Intel\textsuperscript{\textregistered\ } Xeon Silver 4309Y CPU and 64~GB~of~memory.
\subsection{Initial BT MOR}
\label{sec:exp_bt_init}
\label{sec:small_scale_models}
For the evaluation of the initial BT MOR method, we used small-scale RC and RLCk models (i.e., real transmission lines) extracted by ANSYS RaptorX\texttrademark~\cite{raptX}, which are presented in Table~\ref{table:ANSYS_models}.
\begin{table}[H]
\footnotesize 
\centering
\vspace{-0.2em}
	\caption{ Small-scale RC and RLCk models of transmission lines}
        \vspace{-6pt}	
 \label{table:ANSYS_models}
\setlength{\tabcolsep}{4pt}
\begin{tabular}{|c|c|c|c|c|c|c|c|}
\hline
 Model & Initial order & \#nodes &\#ports & \#resistors & \#capacitors & \#inductors & \#mutual inductors \\ \hline
 RC\_1& 48 & 48 & 2 & 202 & 273 & 0 & 0 \\ \hline
 RC\_2& 526 & 526 & 6 & 6667 & 6872 & 0 & 0 \\ \hline
 RLCk\_1& 5431 & 3084 & 2 & 2998 & 1282 & 2347 & 136271 \\ \hline
 RLCk\_2& 21800 & 12166 & 6 & 34635 & 31131 & 9634 & 23639237  \\ \hline
\end{tabular}
\vspace{-0.2em}
\end{table}
\label{sec:results_bt_init} 
As can be seen in Table~\ref{table:vanilla_results}, the mean relative error (MRE) and maximum relative error (MAX\_RE) for the DC analysis are lower than 0.61\% and 0.62\%, respectively. As for the SP analysis, MRE remains under 0.73\% while  MAX\_RE is below 1.89\%. The reduction percentage to achieve accurate results is within 65-87\%. However, considering the runtime and memory overhead of the initial method, its application on large-scale~models is practically~infeasible.
\begin{table}[!hbt]
\footnotesize 
\centering
\vspace{-0.2em}
    \caption{Evaluation of ROMs generated by the initial BT MOR against the original models
    }
      \vspace{-6pt}
    \label{table:vanilla_results}
\setlength{\tabcolsep}{4pt}
\resizebox{\textwidth}{!}{
\begin{tabular}{|c||c|c||c|c|c|c||c||c|}
\hline
\multirow{2}{*}{Model} & \multirow{2}{*}{\makecell{ROM \\ order}} & \multirow{2}{*}{\makecell{Reduction \\ (\%)}} & \multicolumn{2}{c|}{DC analysis} & \multicolumn{2}{c||}{SP analysis} & \multirow{2}{*}{\makecell{Reduction \\ time}} & \multirow{2}{*}{\makecell{Memory \\ (GB)}} \\ 
\cline{4-7} & & & MRE (\%) & MAX\_RE (\%) & MRE (\%) & MAX\_RE (\%) & & \\ \hline 
RC\_1 & 17 & 64.58 & 0.54 & 0.55 & 0.73 & 0.78 & 0.04 s & 0.01 \\ \hline
RC\_2 & 93 & 82.32 & 0.09 & 0.17 & 0.32 & 0.43 & 2.51 s & 0.09 \\ \hline
RLCk\_1 & 1131 & 79.18 & 0.61 & 0.62 & 0.25 & 1.89 & 1.12 h. & 7.24 \\ \hline
RLCk\_2 & 2797 & 87.17 & 0.46 & 0.61 & 0.08 & 1.38 & 6 days & 87.11 \\ \hline
\end{tabular}
}
\vspace{-0.8em}
\end{table}

\subsection{Low-rank BT MOR}
\label{sec:exp_bt_low-rank}
\label{sec:large_scale_models}
To validate the accuracy and performance of the low-rank BT MOR method,
we designed~disparate circuits in the GlobalFoundries 22 nm FDSOI technology and extracted the~corresponding large-scale RLCk models using 
RaptorX\texttrademark~\cite{raptX}. The choice of the benchmark circuits is driven by their diversity and therefore, different metrics are used to describe their behavior. As shown in Table~\ref{table:benchmarks}, the evaluated designs include Hybrid and Wilkinson couplers as well as typical transceiver blocks like low-noise-amplifiers (LNAs) and oscillators, where the metrics of interest are the reflection coefficients and performance (gain, noise, linearity). 
In our experiments, we utilized the MORCIC tool to generate ROMs with target accuracy comparable~to~RaptorX\texttrademark.

\begin{table}[!hbt]
\footnotesize 
\caption{Large-scale RLCk models and metrics of interest for the designed circuits}
\label{table:benchmarks}
\vspace{-6pt}
\resizebox{\textwidth}{!}{
\begin{tabular}{lccll}
\hline
\textbf{Block/DUT}                                                                   & \multicolumn{1}{l}{\textbf{\begin{tabular}[c]{@{}l@{}}Initial\\ Order\end{tabular}}} & \multicolumn{1}{l}{\textbf{Ports}}                & \textbf{\begin{tabular}[c]{@{}l@{}}Mutual \\ inductors\end{tabular}} & \textbf{\begin{tabular}[c]{@{}l@{}}Simulated\\ Metrics\end{tabular}}                                  \\ \hline
\textit{\begin{tabular}[c]{@{}l@{}}Hybrid Coupler @28GHz\\ Hybrid Coupler @56GHz\end{tabular}}    & \begin{tabular}[c]{@{}c@{}} 134710\\ 98024\end{tabular}                                & \begin{tabular}[c]{@{}c@{}}5\\ 5\end{tabular} & \begin{tabular}[c]{@{}l@{}}  79001243\\ 52363149\end{tabular}                        & \begin{tabular}[c]{@{}l@{}} S-parameters of coupler as:\\ power splitter \&  divider\end{tabular} \\ \hline
\textit{\begin{tabular}[c]{@{}l@{}}Wilkinson Coupler @28GHz\\ Wilkinson Coupler @56GHz\end{tabular}} & \begin{tabular}[c]{@{}c@{}} 129087\\ 100888\end{tabular}                                & \begin{tabular}[c]{@{}c@{}} 4\\ 4\end{tabular} & \begin{tabular}[c]{@{}l@{}} 259462454\\ 193641938 \end{tabular}                        & \begin{tabular}[c]{@{}l@{}} S-parameters of coupler as:\\ power splitter \& divider\end{tabular} \\ \hline
\textit{\begin{tabular}[c]{@{}l@{}}VGA @28GHz\end{tabular}}                        & 95189                                                                                & 13                                               & 40230583                                                                    & \begin{tabular}[c]{@{}l@{}}S-parameters,  attenuation\end{tabular}       \\ \hline
\textit{\begin{tabular}[c]{@{}l@{}}VCO @13GHz\end{tabular}}                        & 104367                                                                                & 4                                               & 70445484                                                                    & \begin{tabular}[c]{@{}l@{}}Spectrum, PN, osc. frequency\end{tabular}       \\ \hline
\textit{\begin{tabular}[c]{@{}l@{}}LNA Common-Source @56GHz \\LNA Cascode @28GHz \end{tabular}}          & \begin{tabular}[c]{@{}l@{}}  128574  \\ 162881 \end{tabular}       
&\begin{tabular}[c]{@{}c@{}}  9 \\ 11    \end{tabular}          &\begin{tabular}[c]{@{}c@{}}  72832315    \\ 98585323   \end{tabular}                                                                 & \begin{tabular}[c]{@{}l@{}}S-parameters, gain, CP1dB,\\  IIP3, Noise Figure (NF) \end{tabular}                 \\ \hline

\end{tabular}
}
\end{table}

\label{sec:results_bt_low-rank}

\begin{center}
\begin{figure}[!hbt]
  \vspace{-0.5em}
  \centering
  \includegraphics[width=0.98\textwidth]{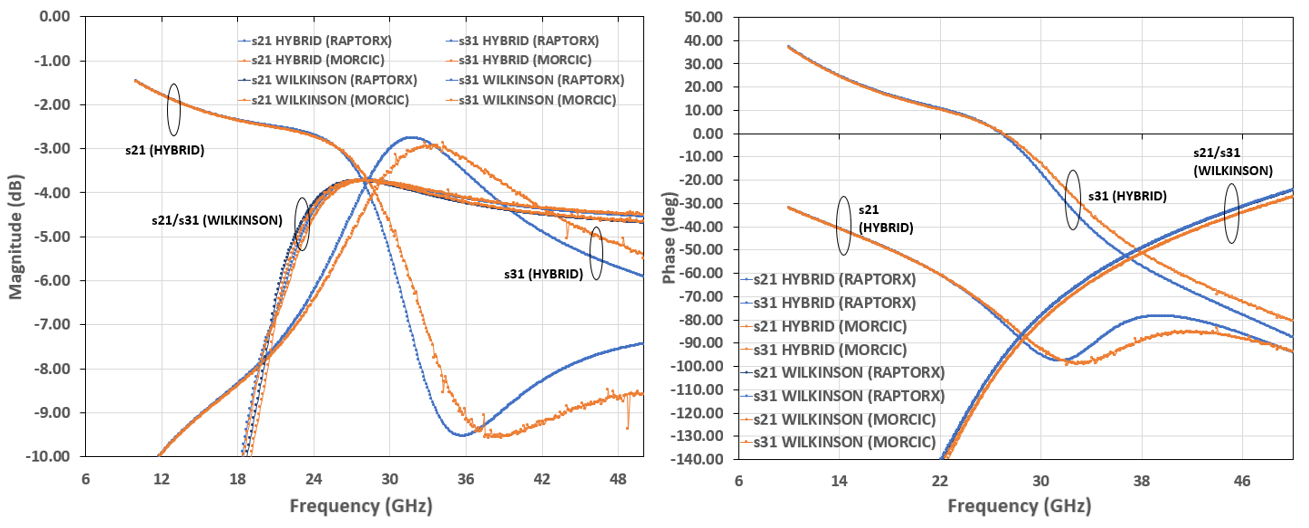}
  \vspace{-7pt}
  \caption{
   MORCIC vs RaptorX\texttrademark\ ROM accuracy: Hybrid/Wilkinson S-parameters and phases.
  }
  \label{fig:morcic_result}
  \vspace{-0.5em}
\end{figure}
\end{center}

The accuracy evaluation 
is performed by comparing the ROMs generated by the low-rank BT MOR against the 
reference ROMs obtained by 
RaptorX\texttrademark, as the simulation of the original extracted models (i.e., full RLCk netlists) is infeasible.
The 
evaluated metrics for the Hybrid~and Wilkinson couplers at 28 GHz, both operating as power splitters, are demonstrated in 
Figure~\ref{fig:morcic_result}. As can be seen, the S-parameters of the MORCIC ROMs closely match 
those 
of the 
RaptorX\texttrademark\ ROMs across the frequency range, and most importantly at the
frequency of interest.
The insertion-loss error is lower than 
0.5 dB, while the respective phases differ by less than 2 degrees.

The efficiency of the MORCIC tool against RaptorX\texttrademark\ is demonstrated in Table \ref{table:low_rank_results}. 
On average, 
MORCIC produces 
$\times$3.1 more compact ROMs.
As the number of ports increases, the advantage of BT is more evident, 
with the maximum improvement in ROM order reaching~$\times$5.5 for 
LNACasc\_28.
Although MORCIC has higher reduction time and memory requirements 
compared to RaptorX\texttrademark,
they are still reasonable and can be significantly improved in future~work.
\begin{table}[!hbt]
\footnotesize 
\centering
\vspace{-1em}
    \caption{
    MORCIC vs RaptorX\texttrademark\ ROM order and MOR performance
}
      \vspace{-6pt}
    \label{table:low_rank_results}
\setlength{\tabcolsep}{5pt}
\begin{tabular}{|c||c||c|c||c|c||c|c|}
\hline
\multirow{2}{*}{Model}  &  \multirow{2}{*}{\makecell{Initial \\ order}} & \multicolumn{2}{c||}{ROM order} & \multicolumn{2}{c||}{Reduction time (s)} & \multicolumn{2}{c|}{Memory (GB)} \\ 
\cline{3-8}& & RaptorX\texttrademark& MORCIC & RaptorX\texttrademark& MORCIC & RaptorX\texttrademark& MORCIC\\ \hline 
VGA\_28 & 95189 & 4744 & 1040 & 67 & 1037 & 32.63 & 19.14 \\ \hline
Hybrid\_56 & 98024 & 1267 & 397 & 104 & 613 & 24.05 & 29.11 \\ \hline
Wilkinson\_56 & 100888 & 765 & 320 & 154 & 570 & 24.79 & 29.76 \\ \hline
VCO\_13 & 104367 & 407 & 311 & 119 & 673 & 26.48 & 29.18 \\ \hline
LNACS\_56 & 128574 & 2172 & 716 & 74 & 1237 & 25.82 & 26.74 \\ \hline
Wilkinson\_28 & 129087 & 885 & 302 & 205 & 801 & 25.35 & 36.21 \\ \hline
Hybrid\_28 & 134710 & 787 & 399 & 217 & 1032 & 24.31 & 35.52 \\ \hline
LNACasc\_28 & 162881 & 4768 & 879 & 373 & 2866 & 78.52 & 48.67 \\ \hline
\end{tabular}
\vspace{-1.4em}
\end{table}
\section{Conclusions}
\label{sec:conclusions}
\vspace{-0.4em}
In this paper, we present efficient BT MOR techniques
to reduce electromagnetic RLCk models. The accuracy of the proposed methods has been evaluated across diverse benchmark circuits, such as the Hybrid/Wilkinson couplers, primarily comparing their S-parameters. Experimental results demonstrate that our low-rank BT MOR approach achieves sufficient accuracy while providing ROMs that are up to ×5.5 smaller than the ROMs obtained by~ANSYS~RaptorX\texttrademark.
\vspace{-0.6em}
\section{Acknowledgments}
\label{sec:acknowledgments}
\vspace{-0.4em}
This research has been co-financed by the European Regional Development Fund and Greek national funds via the Operational Program "Competitiveness, Entrepreneurship and Innovation,"~under the call "RESEARCH-CREATE-INNOVATE" (project code: T2EDK-00609).
\vspace{-1.6em}
\bibliographystyle{IEEEtran}
\bibliography{easychair}

\begin{thebibliography}{1}
\providecommand{\url}[1]{#1}
\csname url@samestyle\endcsname
\providecommand{\newblock}{\relax}
\providecommand{\bibinfo}[2]{#2}
\providecommand{\BIBentrySTDinterwordspacing}{\spaceskip=0pt\relax}
\providecommand{\BIBentryALTinterwordstretchfactor}{4}
\providecommand{\BIBentryALTinterwordspacing}{\spaceskip=\fontdimen2\font plus
\BIBentryALTinterwordstretchfactor\fontdimen3\font minus \fontdimen4\font\relax}
\providecommand{\BIBforeignlanguage}[2]{{%
\expandafter\ifx\csname l@#1\endcsname\relax
\typeout{** WARNING: IEEEtran.bst: No hyphenation pattern has been}%
\typeout{** loaded for the language `#1'. Using the pattern for}%
\typeout{** the default language instead.}%
\else
\language=\csname l@#1\endcsname
\fi
#2}}
\providecommand{\BIBdecl}{\relax}
\BIBdecl

\bibitem{Odabasioglu1998}
A.~Odabasioglu \emph{et~al.}, ``Prima: Passive reduced-order interconnect macromodeling algorithm,'' \emph{IEEE Trans. on CAD of Integrated Circuits and Systems}, vol.~17, no.~8, pp. 645--654, 1998.

\bibitem{Antoulas2004}
S.~Gugercin \emph{et~al.}, ``A survey of model reduction by balanced truncation and some new results,'' \emph{International Journal of Control}, vol.~77, no.~8, pp. 748--766, 2004.

\bibitem{raptX}
\BIBentryALTinterwordspacing
``{Ansys-RaptorX}.'' [Online]. Available: \url{www.ansys.com/products/semiconductors/ansys-raptorh}
\BIBentrySTDinterwordspacing

\bibitem{MNA}
C.-W. Ho \emph{et~al.}, ``{The modified nodal approach to network analysis},'' \emph{IEEE Trans. on Circuits and Systems}, vol.~22, no.~6, pp. 504 -- 509, 1975.

\bibitem{Plancherel}
K.~Gr{\"{o}}chenig, \emph{Foundations of Time-Frequency Analysis}, ser. Applied and numerical harmonic analysis.\hskip 1em plus 0.5em minus 0.4em\relax Birkh{\"{a}}user, 2001.

\bibitem{Lathauwer2004}
D.~Lathauwer \emph{et~al.}, ``Computation of the canonical decomposition by means of a simultaneous generalized schur decomposition,'' \emph{SIAM Journal on Matrix Analysis and Applications}, vol.~26, no.~2, pp. 295--327, 2004.

\bibitem{Simoncini2007}
V.~Simoncini, ``A new iterative method for solving large-scale lyapunov matrix equations,'' \emph{SIAM Journal on Scientific Computing}, vol.~29, no.~3, pp. 1268--1288, 2007.

\bibitem{ASPDAC21}
C.~Chatzigeorgiou \emph{et~al.}, ``{Exploiting Extended Krylov Subspace for the Reduction of Regular and Singular Circuit Models},'' in \emph{Proc. of the 26th Asia South Pacific Design Automation Conference}, pp. 773--778, 2021.

\end{thebibliography}

\end{document}